\documentclass[graybox]{svmult}
\pdfoutput=1

\usepackage{type1cm}        
\usepackage{makeidx}         
\usepackage{graphicx}        
\usepackage{multicol}        
\usepackage[bottom]{footmisc}
\usepackage{newtxtext}       %
\usepackage[varvw]{newtxmath}
\usepackage{csquotes}
\usepackage{amsmath,amsfonts}%
\usepackage{amsthm}%
\usepackage{mathrsfs}%
\usepackage{mathtools}
\usepackage{hyperref}       
\usepackage{cleveref}
\usepackage{titlesec}
\titleformat{\paragraph}[runin] 
  {\normalfont\normalsize\bfseries}{}{0em}{}
\titlespacing{\paragraph}
  {\parindent}{0pt}{0pt}
\crefname{equation}{Eq.}{Eqs.}

\usepackage{enumitem}
\setlist[enumerate]{before=\vspace{0em}}

\usepackage{caption}
\usepackage{tikz}
\DeclareCaptionLabelSeparator{space}{ }
\makeindex             

\begin{document}

\title*{Efficient Numerical Wave Propagation Enhanced By An End-to-End Deep Learning Model}
\titlerunning{Efficient Deep Numerical Wave Propagation}
\author{Luis Kaiser\orcidID{0009-0006-4576-6650} and\\Richard Tsai\orcidID{0000-0001-8441-3678} and\\Christian Klingenberg\orcidID{0000-0003-2033-8204}}
\authorrunning{L. Kaiser et al.}
\institute{Luis Kaiser \at Oden Institute for Computational Engineering and Science, University of Texas at Austin;\\201 E 24th St, Austin, TX 78712, \email{lkaiser@utexas.edu}
\and Richard Tsai \at Oden Institute for Computational Engineering and Science, University of Texas at Austin;\\201 E 24th St, Austin, TX 78712, \email{ytsai@math.utexas.edu}
\and Christian Klingenberg \at Department of Mathematics, University of Wuerzburg; Emil-Fischer-Straße 40, 97074 Wuerzburg, Germany, \email{klingen@mathematik.uni-wuerzburg.de}}

\maketitle

\abstract*{
In a variety of scientific and engineering domains,
the need for high-fidelity and efficient solutions for high-frequency
wave propagation holds great significance.
Recent advances in wave modeling use sufficiently accurate fine solver outputs
to train a neural network that enhances the accuracy of a fast but inaccurate coarse solver.
In this paper we build upon the work of Nguyen and Tsai (2023)
and present a novel unified system that integrates a numerical solver with a deep learning component into an end-to-end framework.
In the proposed setting, we investigate refinements to the network architecture and
data generation algorithm.
A stable and fast solver further allows the use of Parareal, a parallel-in-time algorithm
to correct high-frequency wave components.
Our results show that the cohesive structure improves performance without sacrificing speed,
and demonstrate the importance of temporal dynamics, as well as Parareal, for accurate wave propagation.}

\abstract{In a variety of scientific and engineering domains,
the need for high-fidelity and efficient solutions for high-frequency
wave propagation holds great significance. Recent advances in wave modeling use sufficiently accurate fine solver outputs
to train a neural network that enhances the accuracy of a fast but inaccurate coarse solver.
In this paper we build upon the work of Nguyen and Tsai (2023)
and present a novel unified system that integrates a numerical solver with a deep learning component into an end-to-end framework.
In the proposed setting, we investigate refinements to the network architecture and
data generation algorithm.
A stable and fast solver further allows the use of Parareal, a parallel-in-time algorithm
to correct high-frequency wave components.
Our results show that the cohesive structure improves performance without sacrificing speed,
and demonstrate the importance of temporal dynamics, as well as Parareal, for accurate wave propagation.}

\section{Introduction}\label{ch:introduction}

Wave propagation computations form the forward part of a numerical method for
solving the inverse problem of geophysical inversion.
This involves solving the wave equation with highly varying sound speed many times in a most efficient way.
For instance, by accurately analyzing the reflections and transmissions generated by complex media discontinuities,
it becomes possible to characterize underground formations when searching for natural gas.
However, traditional numerical computations often demand a computationally expensive fine grid
to guarantee stability.

Aside from physics-informed neural networks (PINNs)~\cite{moseley2020solving,MENG2020113250}
and neural operators~\cite{kovachki2023neural,li2020fourier},
convolutional neural network (CNN) approaches yield remarkable results~\cite{nguyen2022numerical,rizzhelm,mlfluid}
to improve the efficiency of wave simulations,
but demand preceding media analysis and tuning of inputs.
Furthermore, numerical solvers are avoided to prioritize speed~\cite{RAISSI2019686};
especially for extended periods, these methods can diverge.

Therefore, combining a classical numerical solver with a neural network to solve the second-order linear wave equation efficiently
across a variety of wave speed profiles is a central point of our research.
We take a first step by expanding the method of Nguyen and Tsai~\cite{nguyen2022numerical} and
build an end-to-end model that enhances a fast numerical solver through deep learning.
Thus, component interplay is optimized,
and training methods can involve multiple steps to account for temporal wave dynamics.
Similarly, while other Parareal-based datasets~\cite{nguyen2022numerical,ibrahim2023parareal}
are limited to single time-steps to add back missing high-frequency components,
a cohesive system can handle Parareal for sequential time intervals.

\paragraph{Approach and Contribution.}~
An efficient numerical solver $\mathcal{G}_{\Delta t} \mathfrak{u} \equiv \mathcal{G}_{\Delta t}[\mathfrak{u},c]$ is used
to propagate a wave $\mathfrak{u}(x,t) = (u, \partial_t u)$ for a time step $t+\Delta t$ on a medium described by the piecewise smooth wave speed $c(x)$ for $x \in [-1,1]^2$.
This method is computationally cheap since the advancements are computed on a coarse grid
using a large time step within the limitation of numerical stability;
however, it is consistently less accurate than an expensive fine solver $\mathcal{F}_{\Delta t} \mathfrak{u} \equiv \mathcal{F}_{\Delta t}[\mathfrak{u},c]$.
Consequently, the solutions $\mathcal{G}_{\Delta t} \mathfrak{u}$ exhibit numerical dispersion errors and miss high-fidelity details.
In a supervised learning framework, we aim to reduce this discrepancy using the outputs from $\mathcal{F}_{\Delta t}$ as the examples.

We define a restriction operator $\mathcal{R}$ which transforms functions from a fine grid to a coarse grid.
Additionally, for mapping coarse grid functions to a fine grid,
we integrate a neural network $\mathcal{I}^\theta$ to augment the under-resolved wave field.
We can now define a neural propagator
$\Psi_{\Delta t}[\mathfrak{u},c, \theta] \equiv \Psi_{\Delta t}^\theta$
that takes a wave field $\mathfrak{u}$ defined on the fine grid, propagates it on a coarser grid,
and returns the enhanced wave field on the fine grid,
\begin{equation}
    \label{eq:eq0}
    \mathfrak{u}_{n+1} \coloneqq \mathfrak{u}(x,t+\Delta t) = \mathcal{F}_{\Delta t} \mathfrak{u}_n \approx \Psi_{\Delta t}^\theta \mathfrak{u}_n \coloneqq \mathcal{I}^\theta \mathcal{G}_{\Delta t} \mathcal{R} \mathfrak{u}_n.
\end{equation}
The models are parameterized by the family of initial wave fields $\mathfrak{F}_{\mathfrak{u}_0}$ and wave speeds $\mathfrak{F}_c$.

\section{Finite-Difference-Based Wave Propagators}\label{ch:theory}

We consider smooth initial conditions and absorbing or periodic boundary conditions
that lead to well-posed initial boundary value problems.
Since we are interested in seismic exploration applications,
both boundary conditions can be used to simulate the propagation of wave fields with initial energy distributed inside a compact domain.
Following the setup in~\cite{nguyen2022numerical}, let $Q_h u$ denote a numerical approximation of $\Delta u$ with discretized spatial and temporal domains, i.e.,
\begin{equation}\label{eq:eq5}
    \partial_{tt} u(x,t) \approx c^2(x) Q_h u(x,t).
\end{equation}
For the spatial ($\Delta x, \delta x$) and temporal spacing ($\Delta t, \delta t$) on uniform Cartesian grids,
the approximation $(u, u_t)_t \approx (u_t, c^2 Q_h u)$ can be solved by a time integrator:
\vspace{-1mm}
\begin{itemize}
\item \textbf{Coarse solver} $\mathcal{G}_{\Delta t^\star}$ $\coloneqq (\mathcal{S}_{\Delta x,\Delta t}^{Q_h})^M$
        with $\Delta t ^\star = M \Delta t$, which operates on the lower resolution grid, $\Delta x \mathbb{Z}^2~\times~\Delta t \mathbb{Z}^+$.
        $Q_h$ is characterized by the velocity Verlet algorithm with absorbing boundary conditions~\cite{enq}.
\item \textbf{Fine solver} $\mathcal{F}_{\Delta t^\star}$ $\coloneqq (\mathcal{S}_{\delta x,\delta t}^{Q_h})^m$ with
        $\Delta t ^\star = m \delta t$, which operates on the higher resolution grid, $\delta x \mathbb{Z}^2~\times~\delta t \mathbb{Z}^+$, and is sufficiently accurate for the wave speed.
        We shall use the explicit Runge-Kutta of forth-order (RK4) pseudo-spectral method~\cite{runge-kutta}.
        Since this approach is only suitable for PDEs with periodic boundary conditions,
        we first apply $\mathcal{F}_{\Delta t^\star}$ to a larger domain and then crop the result.
\end{itemize}

\begin{figure}[!b]
    \centering
    \includegraphics[scale=.78]{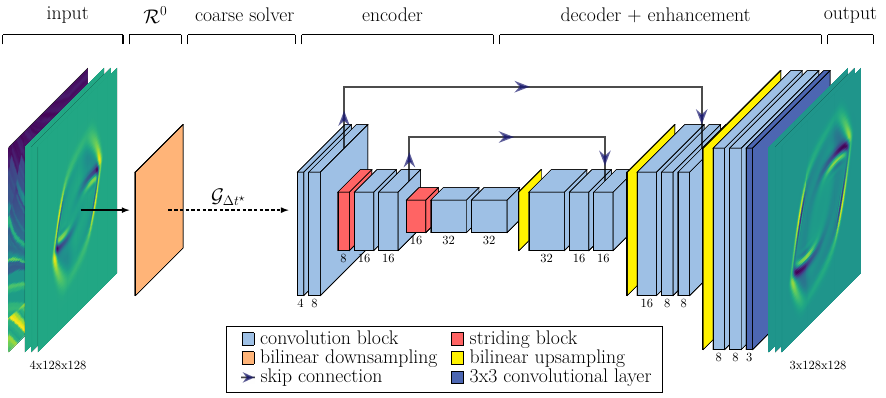}
    \vspace{-5mm}
    \caption{
        Detailed schematic of E2E-JNet3 (adapted from~\cite{nguyen2022numerical}).
        Each convolutional block (blue) encompasses a 3x3 convolutional layer ($\text{groups} = 3$, $\text{padding} = 1$),
        followed by a batch normalization and a ReLU activation function.
        Connectivity within the network is depicted by arrows, with the dashed arrow specifically indicating a single
        application of $\mathcal{G}_{\Delta t^\star}$.
    }
    \label{fig:resnet123}
\end{figure}

\paragraph{Model Components.}~
As the two solvers operate on different Cartesian grids, with $\delta x < \Delta x$ and $\delta t < \Delta t$,
we define the restriction operator $\mathcal{R}$, which transforms functions
from a fine to a coarse grid,
and the prolongation operator $\mathcal{I}$, which maps the inverse relation.
The enhanced variants consist of (a) bilinear interpolations denoted as $\mathcal{R}$ and $\mathcal{I}^{0}$,
while $\mathcal{I}^{0} \mathcal{R} \mathfrak{u} \neq \mathfrak{u}$, and
(b) neural network components
            denoted as $\mathcal{I}^\theta \equiv  \Lambda^{\dag} \mathcal{I}_{\Delta t^\star}^{~\theta} \Lambda$,
            while the lower index indicates that the neural networks are trained when the step size $\Delta t^\star$ is used.
            For improved neural network inference, we use the transition operator
            $\Lambda$ to transform physical wave fields $(u, u_t)$ to energy component representations $(\nabla u, c^{-2}u_t)$, with $\Lambda^{\dag}$ as the pseudo-inverse (see also~\cite{nguyen2022numerical,energy1}).
            \autoref{fig:resnet123} provides a schema visualizing the wave argument transitions.

\paragraph{Variants of the neural propagators.}~
A simple model with bilinear interpolation (E2E-V, $\mathcal{I}^{0} \Psi_{\Delta t^\star}^\theta \mathcal{R}$) is used as a baseline.
Each variant changes the baseline by exactly one aspect.
This allows us to isolate the effect of each architecture modification.
The four investigated end-to-end models $\Psi_{\Delta t^\star}^{\theta} \coloneqq \mathcal{I}_{\Delta t^\star}^{~\theta} \mathcal{G}_{\Delta t^\star} \mathcal{R}$ are:
\begin{description}[leftmargin=!, labelindent=1cm,labelsep=.2cm, labelwidth=2cm, align=left,font=\normalfont\bfseries]
    \vspace{-2mm}
    \item[E2E-JNet3:] E2E 3-level JNet (\autoref{fig:resnet123})
    \item[E2E-JNet5:] E2E 5-level JNet
    \item[E2E-Tira:] Tiramisu JNet~\cite{tiramisu}
    \item[E2E-Trans:] Transformer JNet~\cite{petit2021unet}
    \vspace{-2mm}
\end{description}

The second baseline is taken from~\cite{nguyen2022numerical}
and denoted as the modular, not end-to-end 3-level JNet (NE2E-JNet3),
$\mathcal{I}_{\Delta t^\star}^{~\theta} (\mathcal{G}_{\Delta t^\star} \mathcal{R}$),
while results of $\mathcal{G}_{\Delta t^\star}$ are used to separately train the E2E-JNet3 upsampling component.

\section{Data Generation Approaches}\label{sec:dataset-and-methodology}

For optimal results, the training horizon must be long enough to contain sufficiently representative wave patterns
that develop in the propagation from the chosen distribution of initial wave fields.
Yet the number of iterations must remain small to maintain similarities across different wave speeds.
Similar to~\cite{nguyen2022numerical}, we chose to generate the dataset in the following way:
\begin{enumerate}
  \item An initial wave field $\mathfrak{u}_0 = (u_0, p_0) \in \mathfrak{F}_{u_0}$ is sampled from a Gaussian pulse,
\begin{equation}
    u_0 = e^{- \frac{|x + \tau|^2}{\sigma^2}}, p_0 \equiv \partial_t u_0 = 0
    \label{eq:eq0.1}
\end{equation}
with $x \in [-1,1]^2$, $\frac{1}{\sigma^2} \sim \mathcal{N}(250,10)$ and the initial velocity field $p_0$.
$\tau \in [-0.5,0.5]^2$ is the displacement of the Gaussian pulse's location from the center.
\item Every $\mathfrak{u}_0 \in \mathfrak{F}_{u_0}$ is then propagated eight time steps $\Delta t^\star = 0.06$ by $\mathcal{F}_{\Delta t^\star}$.
  We adopt the fine grid settings for the spatial ($\delta x = \frac{2}{128}$) and temporal resolution ($\delta t = \frac{1}{1280}$) from~\cite{nguyen2022numerical}.
\end{enumerate}

The wave trajectories $\mathfrak{u}_{n+1} = \mathcal{F}_{\Delta t^\star} \mathfrak{u}_n$ provide the input and output data for the supervised learning algorithm,
which aims to learn the solution map $\Psi_{\Delta t^\star}^\theta: X \mapsto Y$:
\begin{equation}
    \begin{split}
        & X \coloneqq \{(\nabla u_n, c^{-2} (u_n)_t, c)\} \\
        & Y \coloneqq \{(\nabla u_{n+1}, c^{-2} (u_{n+1})_t)\},
    \end{split}
    \label{eq:4eq1}
\end{equation}
where $\mathcal{D} = \{(x,y)\}$ with $x \in X$, $y \in Y$.
$\mathcal{D}$ is modified to create $\mathcal{D}^m$, $\mathcal{D}^{w,m}$ (\autoref{subsec:multi-step-training-theory}),
and $\mathcal{D}^p$ (\autoref{sec:parareal-algorithm}).
For brevity, the dataset is only specified if the model is trained on a modified version;
e.g., E2E-JNet3 ($\mathcal{D}^m$) is the E2E-JNet3 model trained on $\mathcal{D}^m$.

\begin{figure}[!t]
    \centering
    \includegraphics[scale=.7, width=\textwidth]{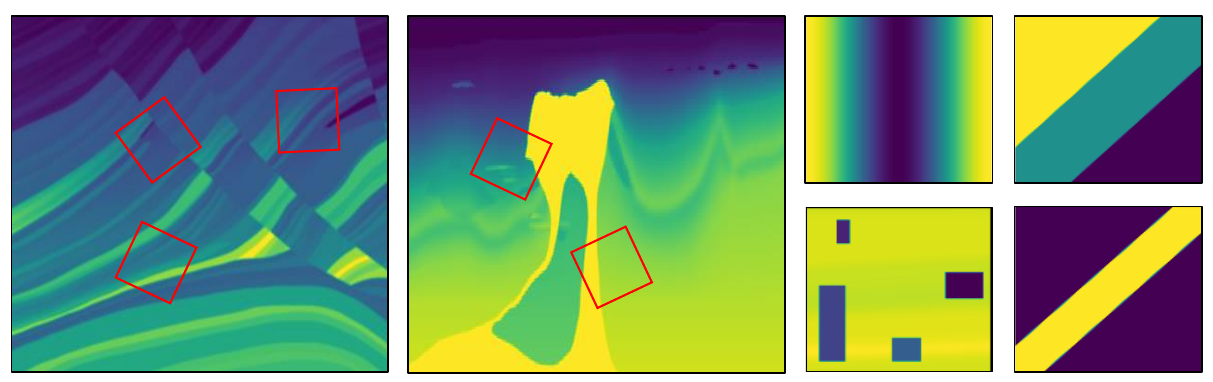}
    \caption{
        Velocity profiles. Brighter colors indicate higher velocity, while randomly chosen subregions are shown in red squares.
        Marmousi and BP profiles are drawn with a probability of $30\%$ each,
        and the other velocity profiles with a probability of $10\%$ each, respectively.
    }
    \label{fig:vel_vis}
\end{figure}

\paragraph{Wave Speeds}~
$c \sim \mathfrak{F}_c$ are sampled from randomly chosen subregions of two synthetic geological structures,
Marmousi~\cite{marmousi} and BP~\cite{bp}, that are mapped onto the spatial grid $h \mathbb{Z}^2 \cap [-1,1]^2$ (see \autoref{fig:vel_vis}).
Four manually modified media (cf.~\cite{nguyen2022numerical}) are added during testing to examine rapid variations in wave speed.


\subsection{Multi-Step Training}\label{subsec:multi-step-training-theory}

During evaluation, the end-to-end model $\Psi_{\Delta t}^\theta$ is applied multiple times to itself to iteratively advance waves over the duration $\Delta t$.
It comes naturally to include longer-term dependencies also in our training dataset.
For $k$ time steps, we therefore introduce a multi-step training strategy that modifies~\cref{eq:eq0}:
\begin{equation}
    \mathfrak{u}_{n+k} \coloneqq \mathfrak{u}(t_n+k\Delta t) = (\mathcal{F}_{\Delta t})^k \mathfrak{u}_n \approx (\Psi_{\Delta t}^\theta)^k \mathfrak{u}_n.
    \label{eq:eq0.0}
\end{equation}
By computing the gradient with respect to the sum of consecutive losses, the gradient flows through the entire computation graph across multiple time steps.

For each initial condition $\mathfrak{u}_0$, $\mathcal{F}_{\Delta t}$ is applied $N$ times
with solutions denoted as $\mathfrak{u}_n, \forall n \in \mathcal{U}_1 \coloneqq \{0,\ldots,N\}$.
In random order, $\Psi_{\Delta t}^\theta$ is applied to every $\mathfrak{u}_n$ for a random amount of steps $k \in \mathcal{U}_2 \coloneqq \{n+1,\ldots,N-n\}$.
Formally, the optimization problem can therefore be described as:
\begin{equation}
    \min_{\theta \in \mathbb{R}^{m \times n}} \mathcal{L}(\Psi_{\Delta t}^\theta;\mathcal{D}) =
    \min_{\theta \in \mathbb{R}^{m \times n}} \frac{1}{\vert \mathcal{D} \vert}
    \sum\limits_{\mathfrak{u}_0,c}
    \sum\limits_{n \sim \mathcal{U}_1 \setminus{\{N\}} } \sum\limits_{\substack{k \sim \mathcal{U}_2  \\
    n < k \leq N-n}} \lVert (\Psi_{\Delta t}^\theta)^k \mathfrak{u}_n - (\mathcal{F}_{\Delta t})^k \mathfrak{u}_n \rVert_{E_h}^2.
    \label{eq:equation2}
\end{equation}
The norm $\lVert \cdot \rVert_{E_h}^2$ is the discretized energy semi-norm MSE as detailed in~\cite{nguyen2022numerical}.
We draw both $n$ and $k$ from the uniform random distributions, i.e., $n \sim \mathcal{U}_1$ and $k \sim \mathcal{U}_2$, respectively.
The novel dataset is denoted as $\mathcal{D}^{m} = \{(\nabla u_n^k, c^{-2} (u_n^k)_t, c; \nabla u_{n+1}, c^{-2} (u_{n+1})_t)\}$.

\paragraph{Weighted Approach.}~
Since in the model's initial, untrained phase, feature variations can be extreme and may lead to imprecise gradient estimations,
we aim to accelerate convergence by weighting individual losses.
Therefore, rather than drawing $k \sim \mathcal{U}_2$ from a uniform distribution, we select values according to a truncated
normal distribution TN$(\mu, \sigma, a, b)$ from the sample space represented as $-\infty < a < b < \infty$.
This focuses on minimizing the impact of errors in the early training stage.
After every third epoch, the mean $\mu$ is increased by one to account for long-term dependencies.
We refer to this dataset as $\mathcal{D}^{w, m}$.

\subsection{Parareal Algorithm}\label{sec:parareal-algorithm}

Identical to~\cite{nguyen2022numerical}, our implemented scheme iteratively refines
the solution using the difference between $\mathcal{F}_{\Delta t}$ and $\mathcal{G}_{\Delta t}$ for each subinterval $\Delta t$.
In particular, missing high-frequency components occur due to the transition to a lower grid resolution,
or a too simple numerical algorithm.
Therefore, a more elaborate model $\Psi_{\Delta t}^\theta$ is required for convergence.
Formally, we rearrange~\cref{eq:eq0} for the time stepping of $\mathcal{F}_{\Delta t}$, and replace $\mathcal{F}_{\Delta t} (\mathcal{I}\mathcal{R}) \mathfrak{u}_n$
by the computationally cheaper strategy $\Psi_{\Delta t}^\theta$ end-to-end:
\begin{gather}
    \mathfrak{u}_{n+1}^{k+1} \coloneqq \Psi_{\Delta t}^\theta \mathfrak{u}_n^{k+1} + [\mathcal{F}_{\Delta t} \mathfrak{u}_n^k - \Psi_{\Delta t}^\theta \mathfrak{u}_n^k],~~~k = 0,\dots,K-1 \label{eq:eq11.24}\\
    \mathfrak{u}_{n+1}^0 \coloneqq \Psi_{\Delta t}^\theta \mathfrak{u}_n^0,~~~n=0,\dots,N-1.\label{eq:eq11.25}
\end{gather}
We observe that the computationally expensive $\mathcal{F}_{\Delta t} \mathfrak{u}_n^k$ on the right-hand side of~\cref{eq:eq11.24} can be performed in parallel for each iteration in $k$.

Parareal iterations alter a given initial sequence of wave fields $\mathfrak{u}_n^0$ to $\mathfrak{u}_n^k$ for $n\in \mathbb{N}$.
This means that neural operators should be trained to map $\mathfrak{u}_n^k$ to $\mathcal{F}_{\Delta t} \mathfrak{u}^k_{n}$.
Therefore, appropriate training patterns for this setup would naturally differ from those found in $\mathcal{D}$,
and the dataset for use with Parareal should be sampled from a different distribution,
denoted as $\mathcal{D}^p$.

\section{Evaluation Setup}\label{ch:evaluation-setup}

The parameters for $\mathcal{G}_{\Delta t^\star}$ are set to $\Delta x = \frac{2}{64}$ and $\Delta t = \frac{1}{600}$,
with a bilinear interpolation scale factor of two.

\paragraph{Experiment 1: Architecture Preselection.}~
The average training time of each variant is approximately $73$ CPU core hours.
Due to resource constraints, we therefore limit our main analysis to one end-to-end variant.
Here, we selected the most promising approach from four deep learning architectures trained on $\mathcal{D}$.

\paragraph{Experiment 2: Multi-Step Training.}~
We train the chosen end-to-end variant from experiment 1 on $\mathcal{D}^{m}$
using an equal number of training points as in $\mathcal{D}$.
The test set is consistent with $\mathcal{D}$ to enable comparison with other experiments.

\paragraph{Experiment 3: Weighted Multi-Step Training.}~
The setup follows experiment 2, while the models are trained on $\mathcal{D}^{w,m}$.

\paragraph{Experiment 4: Parareal Optimization.}~
We explore improvements to our variants using the Parareal scheme in two datasets:\newline
\textit{A. Comprehensive Training} ($\mathcal{D}^p_{\text{train}}$):
The models are trained according to the Parareal scheme in~\cref{eq:eq11.24} and~\cref{eq:eq11.25} with $K=4$
using a random sample that constitutes a quarter of the original dataset $\mathcal{D}$ for fair comparisons.
The gradients are determined by summing the losses of a Parareal iteration.\newline
\textit{B. Fine-tuning} ($\mathcal{D}^p_{\text{refine}}$):
Rather than employing an un-trained model,
we deploy variants that were pre-trained on a random subset containing half of $\mathcal{D}$.
Then, for another subset that constitutes an eighth of $\mathcal{D}$, $\Psi_{\Delta t}^\theta$ is applied with Parareal.

\section{Discussion}\label{ch:discussion}

Each of the total $72$ runs required an average of $72.8$ GPU core hours on one NVIDIA A100 Tensor Core
GPU to complete, while the E2E-JNet3 was trained almost $41\%$ faster and the E2E-Tira three times slower than the average.
This sums up to a total runtime on a single GPU of just over $5{,}241$ hours.

The best trial on the test set was achieved by E2E-Tira with an energy MSE of $0.0109$,
which is well below the $0.0462$ from the previously published model, NE2E-JNet3.
Our most efficient variant is E2E-JNet3 trained on $\mathcal{D}^{w,m}$ with an energy MSE of $0.0169$,
which is close to the results of more extensive models such as E2E-Tira and E2E-Trans,
but is more than five times faster.
\begin{figure}[!b]
    \centering
    \includegraphics[scale=.58]{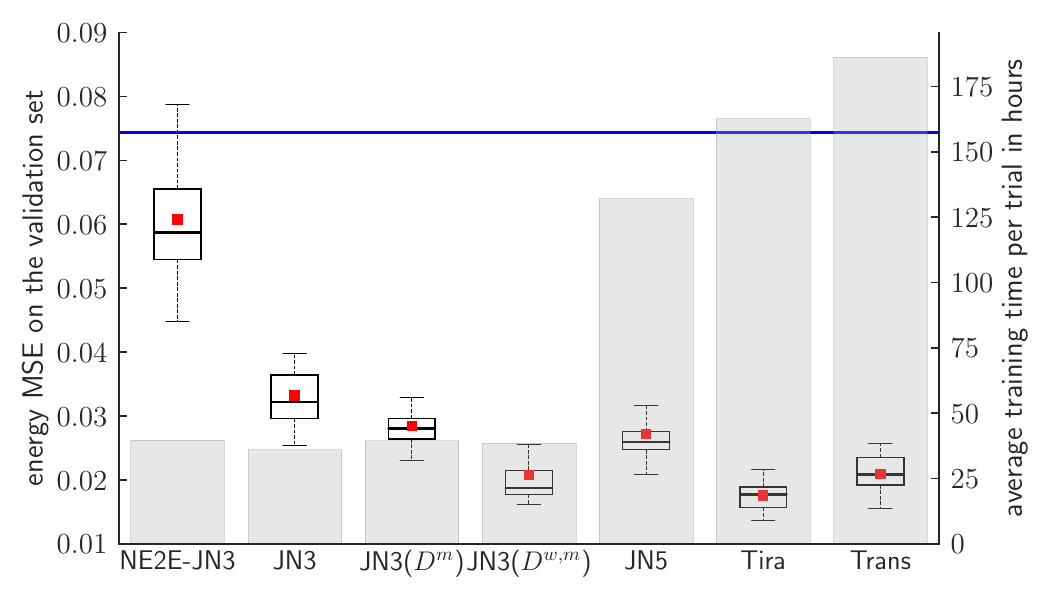}
    \vspace{-2mm}
    \caption{
        Total performance of all hyperparameter search trials on the validation set.
        The boxes represent the range between the 25th and 75th percentile of values,
        while the whiskers indicate 1.5 times the interquartile range.
        The blue line illustrates the result of the baseline E2E-V.
        The red dot shows the mean and the black line marks the median of the data.
        The grey histograms in the background present the average training time of the respective variant in hours on a single GPU.
    }
    \label{fig:box_plot}
\end{figure}

\paragraph{End-to-end Structure.}~
The first important observation based on \autoref{fig:box_plot} is that integrating NE2E-JNet3 into a single, end-to-end system
(E2E-JNet3) improved the average accuracy on the validation set by more than $46\%$, and on the test set by ca. $53\%$.
The ability to include the loss of both the coarse solver and downsampling layer also
caused a lower standard deviation and fewer outliers,
since the mean is well above the median for NE2E-JNet3 compared to E2E-JNet3.

\begin{figure}[!t]
    \centering
    \hspace{-6mm}\includegraphics[scale=.58]{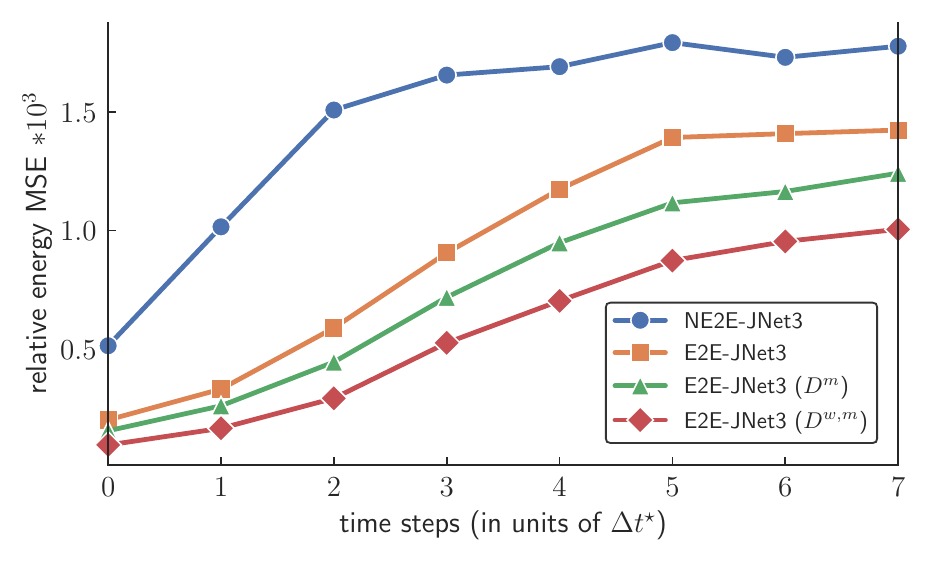}
    \vspace{-2mm}
    \caption{
        Comparing the NE2E-JNet3 model and three end-to-end JNet3 variants that differ in their training algorithms.
        Initial conditions and velocity profiles are sampled from $\mathcal{D}$ and the relative energy MSE results of 10 runs are averaged.
        As expected, all models show a bounded growth as the waves vanish due to absorbing boundary conditions.
    }
    \label{fig:mse_single}
\end{figure}

\begin{figure}[!t]
    \centering
    \includegraphics[scale=.75]{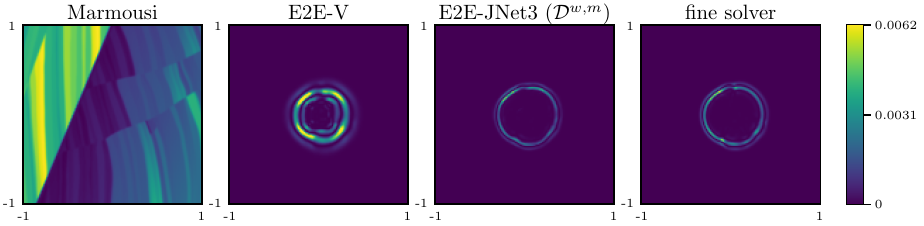}
    \vspace{-5mm}
    \caption{
        Visualizing the wave field correction of E2E-JNet3 ($\mathcal{D}^{w,m}$) in the energy semi-norm after four time steps
        of $\Delta t^\star$.
        The initial condition and velocity profile are sampled from $\mathcal{D}$.
    }
    \label{fig:graphically_compare}
\end{figure}

\paragraph{Multi-Step Training.}~
Introducing a multi-step training loss enhanced
the benefits of an end-to-end architecture even further (cf. E2E-JNet3 ($\mathcal{D}^{m}$) in \autoref{fig:box_plot}) without
increasing the number of model parameters.
\autoref{fig:mse_single} depicts how all end-to-end models had a much lower relative energy MSE (cf.~\cite{nguyen2022numerical})
increase particularly for the first three time steps on the test set.
Hence, we conclude that connecting wave states to incorporate temporal propagation dynamics in the training data appears to be especially important
for the early stages of wave advancements.
Additionally, by taking fewer steps through sampling from a normal distribution that is being shifted along the x-axis (cf. E2E-JNet3 ($\mathcal{D}^{w,m}$)),
we successfully avoid high performance fluctuations when the model is only partially trained.
\autoref{fig:graphically_compare} visualizes the correction of the low-fidelity solution of $\mathcal{G}_{\Delta t^\star}$ by E2E-JNet3 ($\mathcal{D}^{w,m}$).

\paragraph{Upsampling Architecture.}~
An overview of the upsampling architecture performances can be found in~\autoref{fig:table_complexity}.
As expected, the larger networks (E2E-Tira and E2E-Trans) performed slightly better compared to the 3-level JNet architecture,
but for the ResNet architecture (E2E-JNet5), more weights did not increase accuracy by much.
Consequently, we theorize that the ResNet design may be insufficient for capturing high-fidelity wave patterns,
while especially highly-connected layers with an optimized feature and gradient flow (E2E-Tira) are better suited.
Given that E2E-JNet3 ($\mathcal{D}^{w,m}$) had only a slightly worse average energy MSE on the test set,
we generally advise against using the expensive models in our setup.
\begin{table}[!t]
    \centering
    \caption{
        Upsampling variants performance averaged over 10 runs using a batch size of $64$.
    }
    \vspace{-2mm}
    \setlength{\tabcolsep}{0.5em}
\begin{tabular}{l c c c}
 \hline
 \multicolumn{1}{c}{variant} & number of parameters & GPU time (sec) & test energy MSE \\ [0.3ex]
 \hline
 $\mathcal{F}_{\Delta t^\star}$ & - & 57.96749 & - \\
 E2E-V & - & 2.40421 & 0.07437 \\
 E2E-JNet3 & 40{,}008 & 2.88331 & 0.02496 \\
 E2E-JNet5 & 640{,}776 & 10.84893 & 0.02379 \\
 E2E-Tira & 123{,}427 & 13.57449 & 0.01274 \\
 E2E-Trans & 936{,}816 & 15.67633 & 0.01743 \\
 \hline
\end{tabular}

    \label{fig:table_complexity}
\end{table}

\begin{figure}[!t]
    \centering
    \includegraphics[scale=.57]{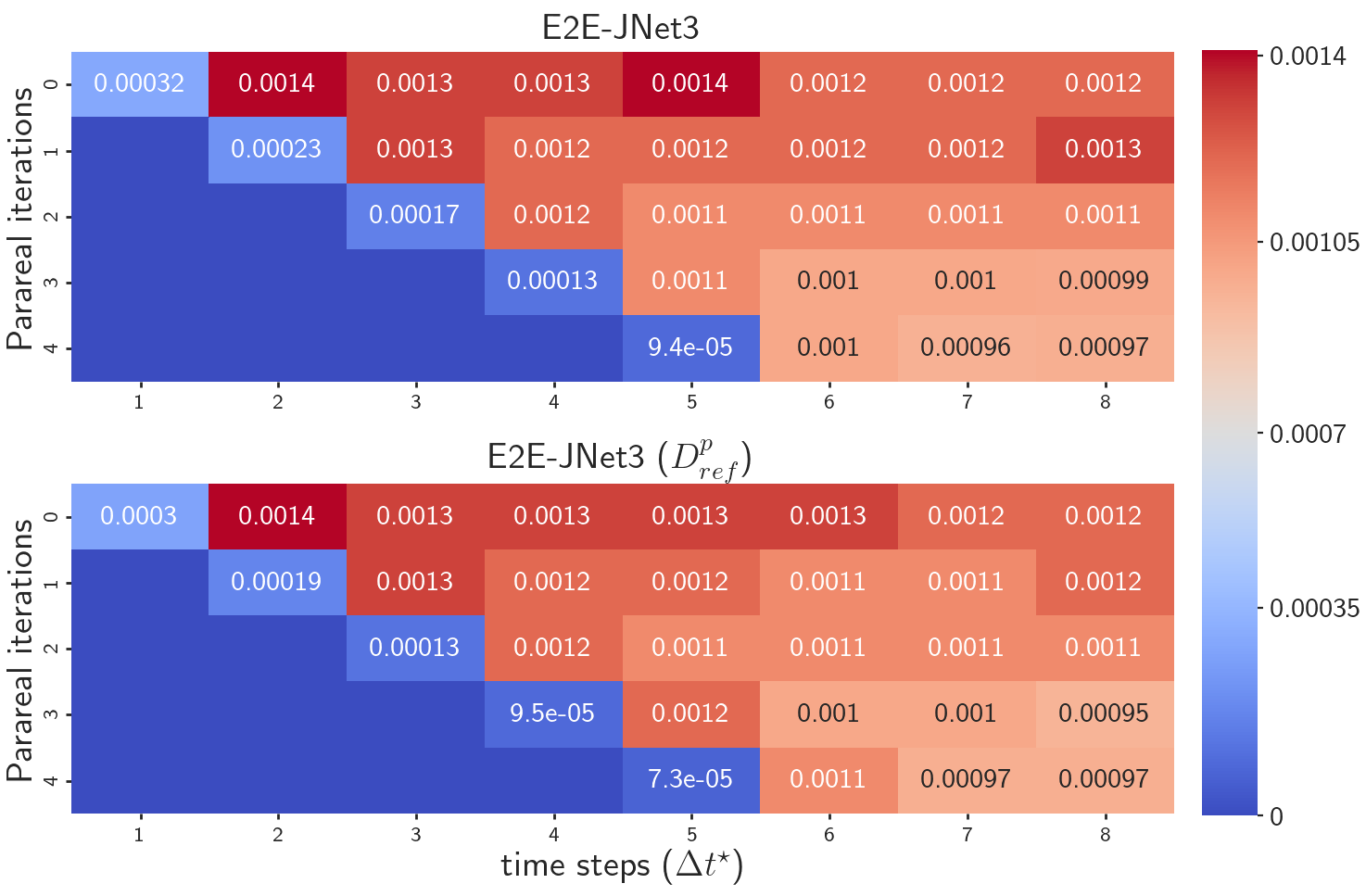}
    \vspace{-1mm}
    \caption{
        Energy MSE of E2E-JNet3 and E2E-JNet3 ($\mathcal{D}^p_{\text{refine}}$) averaged over 10 runs.
    }
    \label{fig:parareal_compare}
\end{figure}

\paragraph{Parareal.}~
While models trained with $\mathcal{D}^p_{\text{train}}$ have an unstable training progress and diverging loss,
applying E2E-JNet3 ($\mathcal{D}^p_{\text{refine}}$) within the Parareal scheme
showed better accuracy than E2E-JNet3 with Parareal (cf.~\autoref{fig:parareal_compare}).
As this training method improved the stability of Parareal,
sampling the causality of concurrently solving multiple time intervals is an efficient enhancement to our end-to-end structure.

\section{Conclusion}\label{ch:conclusion}

In this paper we enhanced the method proposed by Nguyen and Tsai~\cite{nguyen2022numerical},
and reported the results of a large-scale study on different variants that investigate the efficacy of these enhancements.

All end-to-end variants, including the variants with training modifications, outperformed the modular framework of~\cite{nguyen2022numerical}.
In particular, the lightweight end-to-end 3-level JNet (E2E-JNet3) performed reasonably well given its low computation cost,
and was further improved through a weighted, multi-step training scheme ($\mathcal{D}^{w,m}$) to feature time-dependent wave dynamics
without adding complexity to the model or substantially extending the training duration.
Similarly, the Parareal iterations using the neural propagator trained by the Parareal-based data showed significant performance improvements over E2E-JNet3
without extensive additional computational cost due to parallelization.

As expected, certain expensive upsampling architectures, such as intensify the interconnections between feature and gradient flows (Tiramisu JNet),
significantly increased the accuracy.
However, the high computational demand makes its application mostly impractical in modern engineering workflows.


\begin{acknowledgement}
Tsai is partially supported by National Science Foundation Grants DMS-2110895 and DMS-2208504.
The authors also thank Texas Advanced Computing Center (TACC) for the computing resources.
\end{acknowledgement}

\bibliographystyle{spphys}

\begin{thebibliography}{}
\providecommand{\url}[1]{{#1}}
\providecommand{\urlprefix}{URL }
\expandafter\ifx\csname urlstyle\endcsname\relax
  \providecommand{\doi}[1]{DOI \discretionary{}{}{}#1}\else
  \providecommand{\doi}{DOI \discretionary{}{}{}\begingroup
  \urlstyle{rm}\Url}\fi

\end{thebibliography}


\begin{thebibliography}{99.}%
    \bibitem{moseley2020solving} Moseley, B., Markham, A., Nissen-Meyer, T.: Solving the Wave Equation with Physics-Informed Deep Learning. ArXiv eprint:2006.11894 (2020)
    \bibitem{MENG2020113250} Meng, X., Li, Z., Zhang, D., Karniadakis, G. E.: PPINN: Parareal Physics-Informed Neural Network for Time-Dependent PDEs. CMAME \textbf{370}, (2020)
    \bibitem{kovachki2023neural} Kovachki, N., Li, Z., Liu, B., Azizzadenesheli, K., Bhattacharya, K., Stuart, A., Anandkumar, A.: Neural Operator: Learning Maps between Function Spaces with Applications to PDEs. Journal of Machine Learning Research \textbf{24}(89),  pp. 1--97 (2023)
    \bibitem{li2020fourier} Li, Z., Kovachki, N., Azizzadenesheli, K., Liu, B., Bhattacharya, K., Stuart, A., Anandkumar, A.: Fourier Neural Operator for Parametric PDEs. ArXiv preprint:2010.08895 (2020)
    \bibitem{nguyen2022numerical} Nguyen, H., Tsai, R.: Numerical Wave Propagation Aided by Deep Learning. Journal of Computational Physics \textbf{475}, (2023)
    \bibitem{rizzhelm} Rizzuti, G., Siahkoohi, A., Herrmann, F. J: Learned Iterative Solvers for the Helmholtz Equation. 81st EAGE Conference and Exhibition, pp. 1--5 (2019)
    \bibitem{mlfluid} Kochkov, D., Smith, J., Alieva, A., Wang, Q., Brenner, M., Hoyer, S.: Machine Learning–Accelerated Computational Fluid Dynamics.PNAS \textbf{118}(21),  pp. 89--97 (2021)
    \bibitem{RAISSI2019686} Raissi, M., Perdikaris, P., Karniadakis, G.E.: Physics-Informed Neural Networks: A Deep Learning Framework for Solving Forward and Inverse Problems Involving Nonlinear Partial Differential Equations. Journal of Computational Physics \textbf{378},  pp. 686--707 (2019)
    \bibitem{ibrahim2023parareal} Ibrahim, A. Q., Götschel, S., Ruprecht, D.: Parareal with a Physics-Informed Neural Network as Coarse Propagator. ISBN: 978-3-031-39698-4, pp. 649--66 (2023)
    \bibitem{enq} Engquist, B., Majda, A.: Absorbing Boundary Conditions for Numerical Simulation of Waves. PNAS \textbf{74}(5),  pp. 1765--1766 (1977)
    \bibitem{runge-kutta} Runge, C.: Ueber die Numerische Aufloesung von Differentialgleichungen. Mathematische Annalen \textbf{46},  pp. 167–-178 (1895)
    \bibitem{energy1} Rocha, D., Sava, P.: Elastic Least-Squares Reverse Time Migration Using the Energy Norm. Geophysics \textbf{83}(3), pp. 5MJ--Z13 (2018)
    \bibitem{tiramisu} Jégou, S., Drozdzal, M., Vazquez, D., Romero, A., Bengio, Y.: The One Hundred Layers Tiramisu: Fully Convolutional DenseNets for Semantic Segmentation. 2017 IEEE Conference on Computer Vision and Pattern Recognition Workshops (CVPRW), pp. 1175--1183 (2017)
    \bibitem{petit2021unet} Petit, O., Thome, N., Rambour, C., Soler, L.: U-Net Transformer: Self and Cross Attention for Medical Image Segmentation. Machine Learning in Medical Imaging, Springer International Publishing,  pp. 267--276 (2021)
    \bibitem{marmousi} Brougois, A., Bourget, M., Lailly, P., Poulet, M., Ricarte, P., Versteeg, R.: Marmousi, model and data. Conference: EAEG Workshop - Practical Aspects of Seismic Data Inversion (1990)
    \bibitem{bp} Billette, F., Brandsberg-Dahl, S.: The 2004 BP Velocity Benchmark. European Association of Geoscientists \& Engineers (67th EAGE Conference \& Exhibition) (2005)
\end{thebibliography}

\end{document}